\documentclass[12pt, 14paper,reqno]{amsart}
\usepackage{amsmath,amsfonts,amssymb}
\usepackage[breaklinks]{hyperref}
\usepackage{graphicx}
\usepackage{longtable}
\makeatletter
\@namedef{subjclassname@2020}{%
  \textup{2020} Mathematics Subject Classification}
\makeatother

\theoremstyle{plain}
\theoremstyle{plain}

\newtheorem{rem}{Remark}[section]
\numberwithin{equation}{section}
\newtheorem{thm}{Theorem}[section]
\newtheorem{lem}{Lemma}[section]

\newtheorem{prop}{Proposition}[section]

\newtheorem{exa}{Example}[section]

\numberwithin{equation}{section}

\usepackage[breaklinks]{hyperref}
\makeatletter
\let\oldHyPsd@CatcodeWarning\HyPsd@CatcodeWarning
\renewcommand{\HyPsd@CatcodeWarning}[1]{
  \ifnum\pdfstrcmp{#1}{math shift}=0    
  \else                                 
    \oldHyPsd@CatcodeWarning{#1}
  \fi
}
\makeatother

\begin{document} 
\title[Family of elliptic curves]{On the family of elliptic curves $y^2=x^3-m^2x + (pqr)^2$}
\author{Arkabrata Ghosh}
\address{G. Arkabrata@ SRM University-AP, Mangalgiri Mandal, Neerukonda, Amravati, Andhra Pradesh-522502.}
\email{arka2686@gmail.com}

\date{\today}

\keywords{Elliptic curve; rank; Torsion subgroup}
\subjclass[2020] { 11G05, 14G05}

\maketitle

\section*{Abstract}

In this article, we consider a family of elliptic curves defined by $E_{m}: y^2= x^3 -m^2 x + (pqr)^2 $ where $m $ is a positive integer and $p, q, ~\text{and}~ r$ are distinct odd primes and study the torsion as well the rank of $E_{m}(\mathbb{Q})$. More specifically, we proved that if $m \not \equiv 0 \pmod{3}, m \equiv 2 \pmod {2^{k}}$ where $k \geq 5$,  ~\text{and}~ none of the prime numbers $p$, $q$ and $r$ divides $m$,  then the torsion subgroup of $E_{m}(\mathbb{Q})$ is trivial and a lower bound of the $\mathbb{Q}$ rank of this family of elliptic curves is $2$.

\section{Introduction}
\label{sec:1}
The arithmetic of the elliptic curve is one of the most fascinating branches of mathematics as it connects number theory to algebraic geometry. Brown and Myers \cite{BM02} constructed an infinite family $E_m: y^2=x^3-x+m^2$ of elliptic curves over $\mathbb{Q}$ and proved that the family has trivial torsion. Moreover, they showed that $\text{rank}(E_m(\mathbb{Q})) \geq 2$ if $m \geq 2$ and $ \text{rank}(E_m(\mathbb{Q})) \geq 3$ for infinitely many values of $m$. Antoniewiez \cite{An05} considered another family of elliptic curves $C_m: y^2=x^3-m^2x+ 1$ and showed that $\text{rank}(C_m(\mathbb{Q})) \geq 2$ for $m \geq 2$ and $\text{rank}(C_{4k}(\mathbb{Q})) \geq 3$ for the infinite sub-family with $k \geq 1$. Ekinberg \cite{Ei04}, in his Ph.D. thesis, studied the family $D_{m}: y^2=x^3-m^2x + m^2$, and showed that $\text{rank}(D_{m}(\mathbb{Q})) =2$ for all $m \geq 2$. Later Tadic \cite{Ta12} gave a parametrization on a family of the elliptic curve  $C: Y^2= X^3-T^2X +1$ over function fields. He proved that the torsion subgroup of $E$ over the function field $\mathbb{C}(m)$ is trivial and rank of $E$ is $3~\text{and}~4$ over the function fields $\mathbb{Q}(m)$ and $\mathbb{C}(m)$ respectively. Later, Tadic \cite{Tad12} found a family of elliptic curves $E_m:Y^2=x^3-x+T^2$ having $\text{rank}~\geq 3$ over the function field $\mathbb{Q}
(a,i,s,n,k,l)$, where $s^2=i^3 + a^2$. In addition, using the results of \cite{Tad12}, he proved the existence of families with $\text{rank}~\geq 3$ and $ \text{rank}~ \geq 4$ over the fields of rational functions in four variables. Fujita and Nara \cite{FN18} showed that the lower bound of the family of elliptic curves $E_{a,m} ; y^2 = x^3 -a^2x + m^2$ is at least two. Later Juyal and Kumar\cite{JK18} considered the family $E_{m,p}:y^2= x^3-m^2x + p^2$ and showed the lower bound for the rank of $E_{m,p}(\mathbb{Q})$ is $2$.  Finally Chakraborty and Sharma \cite{CS22} considered the family $E_{pq}: y^2=x^3 -m^2 x + (pq)^2$ where $p$ and $q$ are distinct odd primes and showed that it has trivial torsion and rank at least two.

In this article, we consider the family of elliptic curves defined by $E_{m}:  y^2= x^3- m^2x + (pqr)^2$ for distinct odd primes $p$, $q$ and $r$ and certain conditions on the integer $m$. The main aim of this article is to prove the following Theorems.

\begin{thm}
    \label{thm:1.1}
     Let 
  \begin{equation*}
      E_{m}: ~~ y^2= x^3- m^2x + (pqr)^2
  \end{equation*}
  be a family of elliptic curves where $m $ is a positive integer such that $m \not \equiv 0 \pmod{3} ~\text{and}~ m \equiv 2 \pmod {2^{k}}$ where $k \geq 5$, $ p, q, ~\text{and}~ r$ are distinct odd primes and none of them divides $m$. Then $E_{m}(\mathbb{Q})_{tors}= \{\mathcal{O}\}$ where $E_{m}(\mathbb{Q})_{tors}$ denotes the torsion of $E_{m}(\mathbb{Q})$ and $\mathcal{O}$ is the identity element of the group $E_{m}(\mathbb{Q})$.
  
\end{thm}

\begin{thm}
    \label{thm:1.2}
  Let 
  \begin{equation*}
      E_{m}: ~~ y^2= x^3- m^2x + (pqr)^2
  \end{equation*}
  be a family of elliptic curves where $m $ is a positive integer such that $m \not \equiv 0 \pmod{3}  ~\text{and}~ m \equiv 2 \pmod {2^{k}}$ where $k \geq 5$ and $ p, q, ~\text{and}~ r$ are distinct odd primes. Then $\text{rank}~ (E_{m}(\mathbb{Q})) \geq 2$.
\end{thm}

\section{Preliminaries}
\label{sec:2}

In this section, we shall recall some basic facts regarding elliptic curves. Throughout this article, we shall consider the following family of elliptic curves 

\begin{equation}
    \label{eq:2.1}
    y^2=x^3-m^2x + (pqr)^2
\end{equation}

and denotes it by $E_m$.

\subsection{Elliptic curve} 
Let $K$ be a number field. An elliptic curve $E$ over $K$ is an algebraic curve defined by the Weierstrass equation

\begin{equation*}
    E: y^2= x^3 + bx +c ~\text{with}~~~~b,c \in K,
\end{equation*}

and $\Delta=-(4a^3 + 27b^2) \neq 0$. In other words, the above condition is the same as the cubic equation $x^3+bx+c=0$ having three distinct roots in $K$. Moreover, an elliptic curve can be thought of as a smooth (non-singular at every point) algebraic curve of genus one in the projective space $\mathbb{P}^2(K)$, defined by the homogeneous equation $y^2z=x^3 + bxz^2 + cz^3$ and the point $[0:1:0] $ on the curve is denoted by $\mathcal{O}$. Let $E(K)$ denote the set of all $K$-rational points on $E$ with the "point at infinity" $\mathcal{O}$. In other words, 

\begin{equation*}
    E(K)=\{(x,y) \in K \times K: y^2=x^3 + bx+ c \} \cup \{\mathcal{O}\}.
\end{equation*}

Now we can state the following Proposition regarding $E(K)$ which is as follows.

\begin{prop}(\cite{ST92}, Chapter III, section $5$, page $88$)
\label{prop:2.1}
    $E(K)$- The set of $K$ rational points form a finitely generated abelian group, where we denote the group law with $\bigoplus$ (see \cite{ST92}, Chapter I, section $4$, page $30)$ and the additive identity $\mathcal{O}$. 
\end{prop}

The abelian group $E(K)$ is known as the Mordell-Weil group of $E$ over $K$. When $K=\mathbb{Q}$, a result was proved by Mordell regarding $E(\mathbb{Q})$ which Weil then improved for any number field. Together, this result is known as the Mordell-Weil Theorem which states that

\begin{equation*}
    E(K) \cong E(K)_{tors} \oplus \mathbb{Z}^r.
\end{equation*}

Here $E(K)_{tors}$ is the torsion part of $E$ which is a finite abelian group consisting of elements of finite order and the non-negative integer $r$ is called the rank of the elliptic curve which gives us the information about the number of independent points of infinite order of $E$ over $K$.

The structure of a torsion subgroup of an elliptic curve over $\mathbb{Q}$ is well-understood. Mazur \cite{Ma77} and Nagel-Lutz \cite{ST92} Theorems provide a comprehensive knowledge of the torsion subgroup of an elliptic curve over $\mathbb{Q}$. On the other hand, it is challenging to compute the rank of an elliptic curve and as of today, there is no well-defined algorithm to find it. There are plenty of computational ways that have been developed to find the rank but most of them are either computationally complex or involve heavy mathematical machinery. So an efficient way to compute the rank of an elliptic curve is still unknown.

\section{\text{The torsion subgroup of}~ $\texorpdfstring{E_n}{}$}
\label{sec:3}

The principal aim of this section is to prove some results that will be used to prove the Theorem \eqref{thm:1.1} in the section \eqref{sec:4}. So to achieve our goal, we use the technique of reduction modulo a prime which does not divide the discriminant of an elliptic curve say $E$. If any prime  $p'$ which does not divide $\Delta({E})$, then all the roots of the cubic $x^3 + bx +c$ in $\overline{\mathbb{F}}_{p'}$ are distinct and we can say $E$ is an elliptic curve over $\mathbb{F}_{p'}$ or in other words, at $p'$, $E $ has a good reduction. Now given a good reduction of $E$ at $p'$, the application of Theorem \eqref{thm:3.1} gives an injective map from the group of rational torsion points $E(\mathbb{Q})_{tors}$ into the group $E(\mathbb{F}_{p'})$.

\begin{thm}(\cite{Hu87}, Theorem $5.1$) Let $E$  be an elliptic curve over $\mathbb{Q}$. The restriction of reduction homomorphism $r_{p}|_{E(\mathbb{Q})_{tors}}: E(\mathbb{Q})_{tors} \rightarrow E
_{p}( \mathbb{F}_{p})$ is injective for any odd prime $p$ where $E$ has a good reduction and $r_{2}|_{E(\mathbb{Q})_{tors}}: E(\mathbb{Q})_{tors} \rightarrow E_{2}( \mathbb{F}_{2})$ has the kernel at most $\mathbb{Z}/2\mathbb{Z}$ when $E$ has a good reduction at $2$.
    \label{thm:3.1} 
\end{thm}

We need the following remark before going for tools that will help prove the Theorem \eqref{thm:1.2}.

\begin{rem}
    \label{rem:3.1}
    If $P=(x,y)$ is any point on $E_m$, then by the law of addition for doubling a point on an elliptic curve, we denote $2P=(x', y')$ where values of $x' ~\text{and}~ y'$ are as follows,

    \begin{equation}
        \label{eq:3.1}
        \begin{cases}
        x' = \frac{(x^2 +m^2)^2-8x(pqr)^2}{4y^2},\\
        y'= -y -\frac{3x^2-m^2}{2y}(x-x').
        \end{cases}
    \end{equation}
\end{rem}

We will introduce a couple of lemmas and prove it.

\begin{lem}
    \label{lem:3.1}
    There is no point of order $2$ in $E_{m}(\mathbb{Q})$ for every positive integer $m$ and distinct odd primes $p, q ~\text{and}~ r$.
\end{lem}

\begin{proof}
    Suppose that  $E_{m}(\mathbb{Q})$ contains a point of order two, namely $P=(x,y)$. Then $2P= \{\mathcal{O}\} \Longleftrightarrow  P \Longleftrightarrow  -P \Longleftrightarrow y =0 ~\text{and}~ x \neq 0$.  Therefore,

    \begin{equation}
        \label{eq:3.2}
        x^3-m^2x + (pqr)^2=0.
    \end{equation}
  Since the order of $P$ is finite,  $x $ must be an integer by Nagell-Lutz Theorem(\cite{ST92}, page no $56$, chapter $II$, section $5$). Thus, from the equation \eqref{eq:3.2}, we have

  \begin{equation*}
      m^2= x^2 + \frac{(pqr)^2}{x}.
  \end{equation*}
  
Now the above expression implies that \\
$x \in \{\pm 1, \pm p, \pm q, \pm r, \pm p^2, \pm q^2, \pm r^2, \pm pq, \pm pr, , \pm qr, \pm pq^2, \pm p^2q, \pm (pq)^2,  \pm pr^2, \pm p^2r, \pm (pr)^2 \}$\\

 \text{or}~  $ x 
 \in \{ \pm q^2r, \pm qr^2, \pm (qr)^2, \pm pqr, \pm p^2qr, \pm pq^2r, \pm pqr^2, \pm (pq)^2r, \pm p (qr)^2, \pm (pr)^2 q, \pm (pqr)^2 \}$.

Now  $x= \pm 1 \Longleftrightarrow m^2 = 1 \pm (pqr)^2$ which is impossible. Similarly, $x= \pm pqr \Longleftrightarrow m^2 = (pqr)^2 \pm pqr$ which is again impossible. Using easy computations involving divisibility, one can show that other resulting equations in $m$, $p$, $q$, and $r$ do not have solutions. This is a contradiction to our claim and hence $E_{m}(\mathbb{Q})$ has no point of order $2$.

\end{proof}

\begin{lem}
    \label{lem:3.2}
     There is no point of order $3$ in $E_{m}(\mathbb{Q})$ for every positive integer $m \not \equiv 0 \pmod 3 $ and for distinct odd primes $p, q ~\text{and}~ r$.
\end{lem}

\begin{proof}
Assume, On the contrary that, $E_{m}(\mathbb{Q})$ has a point of order $3$ say $P= (x,y)$. Then $3P= \{\mathcal{O}\}$, or equivalently, $2P=-P$ that implies that $x$-coordinate of $(2P)$= $x$- coordinate of $(-P)$, where $P=(x,y)$ and $2P=(x',y')$(values of $x'$ and $y'$ are given by equation \eqref{eq:3.1}). Therefore after putting the values of $x'$ and upon simplification, we have 

\begin{equation}
  \label{eq:3.3}  
  3x^4 -6m^2 x^2+ 12 (pqr)^2x-m^4=0.
\end{equation}

Using Nagel-Lutz Theorem(\cite{ST92}, page no $56$, chapter $II$, section $5$), one obtains that $x$ is an integer. If we consider the equation \eqref{eq:3.3} modulo $3$, we get $m^4 \equiv 0 \pmod 3$ that implies $m \equiv 0 \pmod 3$. This contradicts our hypothesis and hence, we can say that $E_{m}(\mathbb{Q})$ has no point of order $3$.
\end{proof}

Let $P = (x,y) ~\text{and}~2P=(x',y')$ are points of $E_{m}(\mathbb{Q})$. Then double of the point $2P$ is denoted by $4P=(x'', y'')$, where

\begin{equation}
\label{eq:3.4}
    \begin{cases}
        x'' = \frac{(x'^2 + m^2)^2-8x'(pqr)^2}{4y'^2},\\
        y''= -y'-\frac{3x'^2-m^2}{2y'}(x''-x').
    \end{cases}
\end{equation}

\begin{lem}
    \label{lem:3.3}
 $E_{m}(\mathbb{Q})$ has no point of order $5$ for  $m \equiv 2 \pmod 4$, and for distinct odd primes $p, q ~\text{and}~ r$.
\end{lem}

\begin{proof}
    Suppose $E_m$ has a point of order $5$ say $P$. Then $5P= \{\mathcal{O}\} \Longleftrightarrow 4P= -P$. So $x$ coordinate of $4P$ is same as $x$
coordinate of $(-P)$, where $P=(x,y)$ and $4P=(x'',y'')$($x'' $ and $y''$ are given in equation \eqref{eq:3.4}). Upon simplification after putting the value of $x''$ gives

\begin{equation}
    \label{eq:3.5}
    \begin{split}
    & (C^2-8xD^2)^4 + 32m^2 y^4 ( C^2 -8xD^2)^2 -512 D^2 y^6(C^2-8xD^2)+ 256m^4y^8\\
    & =256xy^6[-2y^2-(3x^2-m^2)(\frac{D^2-8xC^2}{4y^2}-x)]^2,
    \end{split}
\end{equation}

where $C=(x^2 +m^2) ~\text{and}~ D= pqr$.

Now if $x$ is even, then by considering the equation \eqref{eq:3.5} modulo $4$, we get $m \equiv 0 \pmod 4$ which is a contradiction to our assumption $ m \not \equiv 0 \pmod 4$. On the other hand, if we assume $x$ is odd, then again reducing the equation \eqref{eq:3.5} modulo $4$, we have $(x^2 + m^2)^8 \equiv 0 \pmod 4$. As $x ~\text{is odd, we have}~ x^2 \equiv 1 \pmod 4$. So we get the following equation

\begin{equation}
\label{eq:3.6}
    (m^2 +1 )^8 \equiv 0 \pmod 4.
\end{equation}

From the equation \eqref{eq:3.6}, we can conclude that $m \equiv {1,3} \pmod 4$, and this contradicts the hypothesis $m \equiv 2 \pmod 4 $. So our assumption is wrong and hence, we can say that  $E_{m}(\mathbb{Q})$ does not contain any point of order $5$. 

\end{proof}

\begin{lem}
    \label{lem:3.4}
    $E_{m}(\mathbb{Q})$ has no point of order $7$ for $m \equiv 2 \pmod 8$, and for distinct odd primes $p, q ~\text{and}~ r$.
\end{lem}

\begin{proof}
    Let $P=(x,y) \in E_{m}(\mathbb{Q})$ is a point of order $7$. Then\\
    $7P=\{\mathcal{O}\} \Longleftrightarrow 6P=-P \Longleftrightarrow x-\text{coordinate of}~ (6P) \Longleftrightarrow x-\text{coordinate of}~ (-P)$. Therefore, after implementing some elementary simplification, we arrive at the following equation.

    \begin{equation}
    \label{eq:3.7}
    \begin{split}
        & 16y^2y'^4[4y'^2 + (3x'^2-m^2)(x''-x')]^2- (C'^2-8D^2x'-4xy'^2)^2[y'^2(C^2-8D^2x)+ y^2(C'^2-8D'^2x)]\\
        & = 4x^2y^2 y'^2(C'^2-8D^2x'-4xy'^2)^2,
        \end{split}
    \end{equation}
where $C=(x^2 +m^2),~ C'=(x'^2 + m^2) ~\text{and}~ D= pqr$. 

   We reduce the equation \eqref{eq:3.7} modulo $4$ to get 

   \begin{equation}
       \label{eq:3.8}
       -C'^4(y'^2 C^2 + y^2 C'^2) \equiv 0 \pmod 4.
   \end{equation}

   Now we consider two cases:\\ 
   \textit{Case-I}: At first, we consider $x \equiv 0 \pmod 2$. Now by putting the values of $C, C' ~\text{and}~ D$ in the  equation \eqref{eq:3.8} and simplifying, we get 
   \begin{equation*}
       m \equiv 0 \pmod 4.
   \end{equation*}
   This contradicts our assumption $m  \equiv 2 \pmod 8 \not \equiv 0 \pmod 4$.

  \textit{Case-II}: Now we take $x \equiv 1 \pmod 2$. As $x$ is odd, we have $x^2 \equiv 1 \pmod 8$. Now reducing the equation \eqref{eq:3.8} modulo $8$, we obtain 

  \begin{equation}
      \label{eq:3.9}
      -C'^4 (y'^2 C^2 + y^2 C'^{2}) \equiv 4xy^2 y'^2C'^4 \pmod {8}. 
  \end{equation}

Furthermore, from the equation \eqref{eq:3.1}, we get
\begin{equation}
\label{eq:3.10}
    \begin{cases}
        x' = \frac{(m^2+1)^2}{4y^2} \pmod 8,\\
        x'^2 = \frac{(1+m^2)^4}{16y^4} \pmod 8,\\
        x'^4= \frac{(1+m^2)^8}{256y^8} \pmod 8,\\
        y'=  -\frac{(3-m^2)(m^4+ 6m^2-4x-3)}{8y^3} \pmod 8 ,\\
        y'^2 = \frac{(3-m^2)^2(m^2 +1)^2}{64y^6} \pmod 8.
        \end{cases}
\end{equation}

Now by substituting the values of $x', x'^2, x'^4, y' ~\text{and}~ y'^2$ in the equation \eqref{eq:3.9}, we get

\begin{equation}
    \label{eq:3.11}
    (1+m^2)^{16}[4(3-m^2)^2 ( 1+m^2)^6 + (1+m^2)^8] \equiv 0 \pmod 8.
\end{equation}
   Now by assumption, we have $ m \equiv 2 \pmod 8$. Using this, we can deduce from the equation \eqref{eq:3.11} that $5 \equiv 0 \pmod 8$ which is impossible. 
\end{proof}

\section{\text{Proof of Theorem} \ref{thm:1.1}}
\label{sec:4}

We will now complete the proof of the Theorem \eqref{thm:1.1} using the results discussed in the section \eqref{sec:3}.

\begin{proof}
    We know the Discriminant of the family of elliptic curves $E_{m}$ given by the equation \eqref{eq:2.1} is 

    \begin{equation*}
        \Delta(E_m)=16[4m^6- 27 (pqr)^4].
    \end{equation*}

We now split the proof into a couple of different cases and they are as follows.

\textit{Case-I} Firstly, we consider $p,q,r \neq 5$. Then we claim that $5$ does not divide $\Delta(E_m)$ and prove it in the following way: Suppose our assumption is wrong. Hence, $4m^6 \equiv 27(pqr)^4 \equiv 2 (pqr)^4 \pmod 5$. Being odd primes different than $5$, by Fermat's little Theorem, we can deduce that $ p^4 \equiv q^4 \equiv r^4 \equiv 1 \pmod 5$. So, $4m^6 \equiv 2 \pmod 5$. On the other hand, for any integer $m$, $ m^6 \equiv 0,1, ~\text{or}~ 4 \pmod 5$ and hence, $4m^6 \equiv 0, 1, 4 \pmod 5$. So our assumption is wrong and so $ 5 \nmid \Delta(E_m)$. Hence we can say that $E_m$ has a good reduction at $5$. Now while we try to reduce $E_{m}$ in $\mathbb{F}_5$, we get two different scenarios.

(a) If $p^2 \equiv 1 \pmod 5$, then $q^2 \equiv 1, 4 \pmod 5$ and $r^2 \equiv 1,4 \pmod 5$.  It implies that $(pqr)^2 \equiv 1, 4 \pmod 5$. 

(i) When $(pqr)^2 \equiv 1 \pmod 5$, then depending upon whether $m^2 \equiv 0,1 ~\text{or}~ 4 \pmod 5$, $E_m$ reduces into the equations $y^2 =x^3 +1$, $y^2 = x^3 -x + 1$ and $y^2= x^3 + 4x +5 $ respectively. So the corresponding size of the cardinality of $E_{m}(\mathbb{F}_5)$ would be $ 6, 8 ~\text{and}~ 8$.

(ii) When $(pqr)^2 \equiv 4 \pmod 5$, then the curve $E_{m}$ reduces to $y^2 = x^3 + 4$, $y^2 =x^3 -x +4$ and $y^2= x^3 +4x + 4$ according as $m^2 \equiv 0,1 ~\text{or}~ 4 \pmod 5$ respectively with the corresponding cardinality of $E_{m}(\mathbb{F}_5)$ being $ 6, 8 ~\text{and}~ 8$.

(b) If $p^2 \equiv 4 \pmod 5$, then depending upon the choices of $q^2, r^2 \pmod 5$, we will have $(pqr)^2 \equiv 1, 4 \pmod 5$. This is completely analogous to the previous one. 

By Theorem \eqref{thm:3.1} and Lagrange's theorem,  we can say that the possible orders of $E_{m}(\mathbb{Q})_{tors}$ are $1,2,3,4,6,8 ~\text{and}~ 9$ only. Now Lemmas \eqref{lem:3.1} and \eqref{lem:3.2} ensure that $E_{m}(\mathbb{Q})$ does not have any points order $ 2 ~\text{and}~ 3$. Thus we can conclude that 

\begin{equation*}
    E_{m}(\mathbb{Q})_{tors} = \{\mathcal{O}\}.
\end{equation*}

\textit{Case-II} Now we consider the scenario when $p=5$ or $q=5$ or $r=5$ (as $ p \neq q \neq r$). At first, let $p=5$. Then the equation \eqref{eq:2.1} of $E_{m}$ can be re-written as

\begin{equation*}
    y^2= x^3- m^2 x + 25q^2 r^2,
\end{equation*}

with 
\begin{equation*}
    \Delta(E_m)=16[4m^6- 3^3 5^4(qr)^4].
\end{equation*}

(a) if both $ q \neq 7 \neq r$, then it would imply $7$ does not divide $\Delta(E_m)$, and thus it is evident $E_m$ has a good reduction at $7$. Now as both $q$ and $r$ are odd, so $q^2 \equiv r^2 \equiv 1,2 ~\text{or}~ 4  \pmod 7$. Hence, based on different choices of $q$ and $r$, $(qr)^2 \equiv 1, 2 ~\text{or}~ 4  \pmod 7$. Hence, there are three cases to consider while we reduce $E_m$ in $\mathbb{F}_7$.

(i) $(qr)^2 \equiv 1 \pmod 7$:

Depending on whether $m^2 \equiv 0, 1, 2 ~\text{or}~ 4  \pmod 7$, the curve $E_m$ reduces to $y^2 =x^3 + 4$, $y^2= x^3 -x +4$, $y^2= x^3 -4x +4$ or $y^2= x^3 -2x +4$ respectively. So the corresponding cardinality of $E_{m}(\mathbb{F}_7)$ would be $3, 10, 10 ~\text{and}~ 10$ respectively.

(ii) $(qr)^2 \equiv 2 \pmod 7$:

In this case, the curve $E_{m}$ reduces to $y^2 = x^3 + 1$, $y^2 =x^3 -x +1$,  $y^2= x^3 -2x + 1$ or $y^2= x^3 -4x +1$ according as $m^2 \equiv 0,1,2 ~\text{or}~ 4 \pmod 7$ with the corresponding cardinality of $E_{m}((\mathbb{F}_7)$ being $ 12, 12, 12 ~\text{and}~ 12$ respectively .

(iii) $(qr)^2 \equiv 4 \pmod 7$:

Depending upon whether $m^2 \equiv 0, 1, 2 ~\text{or}~ 4  \pmod 7$, the curve $E_m$ reduces to $y^2 =x^3 + 2$, $y^2= x^3 -x +2$, $y^2= x^3 -2x +2$ or $y^2= x^3 -4x +2$ respectively. So the corresponding cardinality of $E_{m}(\mathbb{F}_7)$ would be $9, 9, 9 ~\text{and}~ 9$ respectively.

So from the above three scenarios, it is clear that the possible orders of 
    $E_{m}(\mathbb{Q})_{tors}$ are $1,2,3,4,5, 6,9,10 ~\text{and}~ 12$. Now using the results of lemmas \eqref{lem:3.1}, \eqref{lem:3.2} and \eqref{eq:3.3}, we can say that $E_m$ does not have any points of order $2,3 ~\text{and}~5$. Thus in this case too 

   \begin{equation*}
    E_{m}(\mathbb{Q})_{tors} = \{\mathcal{O}\}.
\end{equation*}

(b) If $q=7$, then the equation \eqref{eq:2.1} of $E_{m}$ can be re-written as

\begin{equation*}
    y^2= x^3- m^2 x + 5^2 7^2 r^2,
\end{equation*}

with 
\begin{equation}
\label{eq:3.12}
    \Delta(E_m)=16(4m^6- 3^3 5^4 7^4 r^4)= 3(21m^6 -3^2 5^4 7^4 r^4) + m^6.
\end{equation}

Now as $r$ is an odd prime, we have $r^2 \equiv 1 \pmod 3$, and from our assumption $ m \neq 0 \pmod 3$, we can say $m ^2 \equiv 1 \pmod 3$ which when combined with the equation \eqref{eq:3.12}, gives $3$ does not divide $ \Delta(E_m)$. So $E_{m}$ has a good reduction at $3$. Now while reducing $E_{m}$ to $\mathbb{F}_{3}$, we have two scenarios that are described below:

(a) At first, when we take $ q \not \equiv 0 \pmod 3$, we have $q^2 \equiv 1 \pmod 3$. Hence, based on the choice of $m$, curve $E_
{m}$ reduces to $ y^2 = x^3-x + 1$ with the corresponding cardinality of $E_{m}(\mathbb{F}_3)$ is $7$.

(b) Now when $q =3$, then  then the curve $E_{m}$ reduces to $y^2 = x^3 -x$. Now the cardinatlity of $E_{m}(\mathbb{F}_3)$ is $4$. 

Now from the above two scenarios, we can say the possible orders of 
    $E_{m}(\mathbb{Q})_{tors}$ are $1, 2, 4 ~\text{and}~ 7$.  Now using the results of lemmas \eqref{lem:3.1}, \eqref{lem:3.2}  \eqref{lem:3.3} and \eqref{lem:3.4}, we can say that $E_m$ does not have any points of order $1, 2, 4  ~\text{and}~7$. Thus in this case too 

   \begin{equation*}
    E_{m}(\mathbb{Q})_{tors} = \{\mathcal{O}\}.
    \end{equation*}

So by combining all the different scenarios, we can conclude that when $m$ is an integer such that $m \not \equiv 0 \pmod{3}$ and $ m \equiv 2 \pmod 8$, and $p, q, r$ are distinct odd primes, then 

\begin{equation*}
    E_{m}(\mathbb{Q})_{tors} = \{\mathcal{O}\}.
    \end{equation*}
    
\end{proof}

\section{\text{The rank of}~ $\texorpdfstring{E_m}{}$}
\label{sec:5}

The rank of an elliptic curve is a topic of major importance in number theory and as of today, it is yet to understand fully. In this section, we will show that the family of the concerned elliptic curve $E_m$ has rank at least two by showing the existence of two linearly independent rational points. We will consider two points namely $A_m=(0, pqr)$ and $B_m=(m, pqr)$ on $E_{m}(\mathbb{Q})$ and need to show that they are linearly independent, i.e, there does not exist two non-zero integers $a$ and $b$ such that

\begin{equation}
\label{eq:4.1}
    [a]A_m + [b]B_m= \mathcal{O}
\end{equation}

where $[a]A_m$ denotes $a$- times addition of $A_m$ in $E_{m}(\mathbb{Q})$.

At first, we claim that the point $B_m=(m, pqr) \in E_{m}(\mathbb{Q})$ is a point of infinite order because if not, then the order of $B_m$ must be finite and hence, it must be in the torsion of  $E_{m}(\mathbb{Q})$. But as we have already proved that 

\begin{equation*}
    E_{m}(\mathbb{Q})_{tors} = \{\mathcal{O}\}, 
    \end{equation*}

we can conclude that the order of $B_m=(m,pqr)$ must be infinite. Hence, we can say the rank of $E_m$ is at least one. We need to use the following result to show the rank is at least two.

\begin{thm}(\cite{Cr97}, Section $3.6$, Page $78$)
\label{thm:4.1}
Let $E(\mathbb{Q})$ (respectively $2E(\mathbb{Q}))$ be the group of rational points (respectively, double of rational points) on an Elliptic curve $E$, and suppose $E$ has trivial torsion. Then the quotient group $E(\mathbb{Q})/2E(\mathbb{Q})$ is an elementary abelian $2$- group of order $2^r$, where $r$ is the rank of $E(\mathbb{Q})$.
\end{thm}

Now to accomplish our aim, we need to prove the following four auxiliary results.

\begin{lem}
    \label{lem:4.1}
    Let $A=(x'.y')$ and $B=(x,y)$ be points on $E_{m}(\mathbb{Q})$ such that $A=2B$ and $x' \in \mathbb{Z}$. Then
    \begin{itemize}
        \item $x \in \mathbb{Z}$,
        \item $x \equiv m \pmod 2$.
    \end{itemize}
\end{lem}

\begin{proof}
Substituting $x=\frac{u}{v}, ~\text{with}~ gcd(u,v)=1 $ in the expression of $x'$ of equation \eqref{eq:3.1} and after simplification, we get

\begin{equation}
    \label{eq:4.2}
    (m^4-4p^2q^2 r^2 x') v^4 + (4m^2x'-8p^2q^2r^2)uv^3 + 2m^2 u^2 v^2-4x'u^{3}v+u^4=0. 
\end{equation}
Now from the equation \eqref{eq:4.2}, it is evident that $v|u^4$ and therefore $v=\pm 1$. Thus $x \in \mathbb{Z}$.

Now the equation \eqref{eq:3.1} can be written as

\begin{equation*}
    (x^2 + m^2)^2=4[x'(x^3-m^2x+ (pqr)^2) + 2x(pqr)^2], 
\end{equation*}

which implies that $2|(x^2 + m^2)$. Thus $x \equiv m \pmod 2$.

\end{proof}

\begin{lem}
\label{lem:4.2}
    The equivalence class $[A_m]=[(0,pqr)]$ is a non-zero element of $E_{m}(\mathbb{Q})/ 2 E_{m}(\mathbb{Q})  $ for any positive integer $m \equiv 2 \pmod {2^5}$ and for odd primes $p,q, ~\text{and}~r $.
\end{lem}

\begin{proof}
    Suppose $A_m=2C$ for some $C=(x,y) \in E_{m}(\mathbb{Q})$. From the equation \eqref{eq:3.1}, we have
    \begin{equation*}
        \frac{(x^2+m^2)^2-8x(pqr)^2}{4y^2}=0,
    \end{equation*}
        which upon simplification, becomes
        \begin{equation}
            \label{eq:4.3}
            (x^2+m^2)^2= 8x(pqr)^2.
        \end{equation}

    The left-hand side of the equation \eqref{eq:4.3} is a perfect square and so the right-hand side must be too. This implies that $x=2k_1^2$ for some $k_1 \in \mathbb{Z}$. Now we claim that $gcd(2,k
_1)=1$. We prove our claim in the following way: Suppose our claim is false. Then $k_1=2^{\alpha} k_2 ~\text{where}~ \alpha, k_2 \in \mathbb{Z} ~\text{and}~ \alpha \geq 1$. Now substituting $x=2k_1^2=2^{2\alpha +1}k_2^2$ in the equation \eqref{eq:4.3}, we obtain

\begin{equation*}
    2^{8\alpha +4}k_2^8+ m^4 + 2^{4\alpha +3} k_2^4 m^2= 2^{2\alpha +4}k_2^2 (pqr)^2
\end{equation*}

which can be re-written as

\begin{equation}
    \label{eq:4.4}
    m^4 + 8 \times 2^{4\alpha}k_2^4 m^2  + 16 \times 2^{8\alpha}k_2^4= 16 \times 2^{2\alpha}k_2^2(pqr)^2.
\end{equation}

Now if we consider the equation \eqref{eq:4.4} modulo $32$, we get $ m^ 4 \equiv 0 \pmod {32}$. Now from our assumption, we can say $m \equiv 2 \pmod {32}$ from which we can deduce that $m^4 \equiv 16 \pmod {32}$. Hence we arrive at a contradiction.


Now by putting the value of $x$ in the equation \eqref{eq:4.3} and simplifying, we obtain

\begin{equation}
    \label{eq:4.5}
    16k_1^8 + m^4 +8k_1^4 m^2=16k_1^2 (pqr)^2.
\end{equation}

Now from the hypothesis, we get $m \equiv 2 \pmod{32}$ that implies $m ^2 \equiv 4 \pmod {32}$ and $m^4 \equiv 16 \pmod {32}$. So when we consider the equation \eqref{eq:4.5} modulo $32$, get that
\begin{equation*}
    k_1^ 8 + 1 -k_1^2 (pqr)^2 \equiv 0 \pmod 2.
\end{equation*}
 
As $k_1$ is odd and $p,q `\text{and} r$ are odd primes, we have $k_1^2 \equiv 1 \pmod 2$, $p^2 \equiv 1 \pmod 2$, $q^2 \equiv 1 \pmod 2$ and $r^2 \equiv 1 \pmod 2$. Substituting these values in the above equation, we get a contradiction that $ 1 \equiv 0 \pmod 2$. So we can say that the equation \eqref{eq:4.5} has no solution. Therefore, $A_m \not \in 2 E_{m}(\mathbb{Q})$.

\end{proof}

\begin{lem}
    \label{lem:4.3} The equivalence class $[B_m]=[(m, pqr)]$ is a non-zero element of $E_{m}(\mathbb{Q})/ 2 E_{m}(\mathbb{Q})$ for positive integers $ m \equiv 2 \pmod {2^{5}}$ and for odd primes $p,q ~\text{and}~ r$.
\end{lem}

\begin{proof}
    Assume $B_m=(m,pqr)=2C$ for some $C=(x,y) \in E_{m}(\mathbb{Q})$.  Now from the equation \eqref{eq:3.1}, we can say

    \begin{equation*}
        \frac{(x^2+m^2)^2 -8x(pqr)^2}{4y^2}=m.
    \end{equation*}

Now from the Lemma \eqref{lem:4.1}, we get $x \equiv m \pmod 2$. Now by putting $x-m=2s$ in the above equation and simplifying, we get

\begin{equation}
    \label{eq:4.6}
    (2s^2-m^2)^2=(pqr)^2[4s + 3m].
\end{equation}
    Now as the left-hand side of the equation \eqref{eq:4.6} is a perfect square, so would be the right-hand side. So, we have $(4s+3m)=w^2$ for some $w \in \mathbb{Z}$ which in turn implies $3m \equiv {0,1} \pmod 4$. Now as $ m \equiv 2 \pmod {2^5}$, we have $m \equiv 2 \pmod 4$. So from here, we get $3m \equiv 2 \pmod 4$ which is a contradiction to the fact that $(4s+3m)$ is a perfect square. Hence, $B_m \notin 2 E_{m}(\mathbb{Q})$.
\end{proof}

\begin{lem}
    \label{lem:4.4}
    The equivalence class $[A_m+ B_m]=[-m, pqr]$ is a non-zero element of $E_{m}(\mathbb{Q})/ 2 E_{m}(\mathbb{Q})$ for positive integers $ m \equiv 2 \pmod {2^{5}}$ and for odd primes $p,q ~\text{and}~ r$.
\end{lem}

\begin{proof}
    Suppose $A_m + B_m =(-m, pqr)=2C$ for some $C=(x,y) \in E_{m}(\mathbb{Q})$. Thus,
       \begin{equation*}
           \frac{(x^2+m^2)^2 -8x(pqr)^2}{4y^2}=-m.
       \end{equation*}
  As $x \equiv m \pmod 2$, we can write $x-m=2s$. Now substituting this value in the above equation and  after simplification, we have

  \begin{equation}
  \label{eq:4.7}
      4s^4+ 16ms^3 +20m^2s^2 + 8m^3s -4s(pqr)^2 + m^4 -(pqr)^2m=0
  \end{equation}

Now as $m \equiv 2 \pmod {2^{5}}$, we can say $m \equiv 2 \pmod {16}$. Now as $m \equiv 2 \pmod {16}$, we have $m^2 \equiv 4 \pmod {16}$, $m^{3} \equiv 8 \pmod {16}$ and $m^{4} \equiv 0 \pmod {16}$. Now substituting these values in the equation \eqref{eq:4.7}, we get

\begin{equation*}
    4s^4 -4s(pqr)^2 -2 (pqr)^2  \equiv 0 \pmod {16},
\end{equation*}

which can be simplified into the equation

\begin{equation}
    \label{eq:4.8}
   2 s^4- 2 s (pqr)^2 -(pqr)^2 \equiv 0 \pmod 8.
\end{equation}

Now as $p, q ~\text{and}~ r$ are odd primes, we have $p^2\equiv q^2 \equiv r^2 \equiv 1 \pmod 8$. Then, from equation \eqref{eq:4.7}, we have $2s[s^3-(pqr)^2] \equiv 1 \pmod 8$ which is a contradiction.

\end{proof}

Now if we can show that $\{\mathcal{O}, [A_m], [B_m], [A_m + B_m]\}$ is a subgroup of $E_{m}(\mathbb{Q})/ 2 E_{m}(\mathbb{Q})$ and $A_m, ~B_m$ are linearly independent, then we can prove that rank of $E_{m}(\mathbb{Q})$ is at least $2$. Now to prove the Theorem \eqref{thm:1.2}, we will state and prove the following propositions.

\begin{prop}
    \label{prop:4.1}
    Let $m$ be a positive integer such that $m \not \equiv 0 \pmod{3} ~\text{and}~ m \equiv 2 \pmod {2^{k}}$ where $k \geq 5$ with odd primes $p,q ~\text{and}~ r$. Then the set 
    \begin{equation*}
         S=\{\mathcal{O}, [A_m],[B_m], [A_m + B_m]\}
    \end{equation*}
    is a subgroup of $E_{m}(\mathbb{Q})/ 2 E_{m}(\mathbb{Q})$ of order $4$ with $A_m=(0, pqr),~ B_m=(m, pqr)$.
\end{prop}

\begin{proof}
    From the Lemmas \eqref{lem:4.2}, \eqref{lem:4.3} and \eqref{lem:4.4}, we can say that all of $[A_m], [B_m] ~\text{and}~ [A_m + B_m]$ are not equal to $\{\mathcal{O}\}$. Now if $[A_m]=[B_m]$, then $[A_m + B_m]=[2B_m]=[\{\mathcal{O}\}]$, which is not possible. Now if $[A_m]=[A_m + B_m]$, then we have $[\{\mathcal{O}\}]=[2A_m]=[B_m]$, which is again a contradiction. Similarly, we can show $[B_m]$ and $[A_m + B_m]$ are also distinct. Hence, $\mathcal{O}, [A_m], [B_m] ~\text{and}~ [A_m + B_m] $ are distinct classes in $E_{m}(\mathbb{Q})/ 2 E_{m}(\mathbb{Q})$. Thus the set $S$ is a subgroup of order $4$ in $E_{m}(\mathbb{Q})/ 2 E_{m}(\mathbb{Q})$ .
\end{proof}

\begin{prop}
    \label{prop:4.2}
    Points $A_m=(0,pqr)$ and $B_m=(m, pqr)$ are two linearly independent points in $E_{m}: y^2= x^3-m^2 x +(pqr)^2$ where $ m \neq 0 \pmod 3$ and $ m \equiv 2 \pmod {2^{k}},~ k\geq 5$ and for odd primes $p,q ~\text{and}~ r$. 
\end{prop}

\begin{proof}
    From the discussion at the beginning of this section, we know that if the points $A_m$ and $B_m$ are linearly independent, then there do not exist non-zero integers $a$ and $b$ such that it satisfies equation \eqref{eq:4.1}. Suppose, on the contrary, we have $aA_m + bB_m= \mathcal{O}$ where $a$ and $b$ are defined as above with $a$ is minimal. We need to consider four different cases and they are as follows:

    \begin{itemize}
        \item If $a$ is even and $b$ is odd and $[aA_m + bB_m]= [\mathcal{O}]$, then in the group $ E_{m}(\mathbb{Q})/ 2 E_{m}(\mathbb{Q})$, we have $[B_m]=[\mathcal{O}]$. This is a contradiction by Lemma \eqref{lem:4.3}.

        \item If $a$ is odd and $b$ is even, then $[aA_m + bB_m]= [\mathcal{O}]$ implies that $[A_m]=[\mathcal{O}]$. This is impossible due to the lemma \eqref{lem:4.2}.

        \item When both $a$ and $b$ are odd we get $[A_m + B_m]= [\mathcal{O}]$, which contradicts lemma \eqref{lem:4.3}.

        \item Now if both $a$ and $b$ are even, then we may assume $a=2a'$ and $b=2b'$. Then from $[aA_m + bB_m]= [\mathcal{O}]$, we get that $ 2[a'A_m + b'B_m]= [\mathcal{O}]$, which implies that $(a'A_m + b'B_m)$ is a point of order $2$ in $E_{m}(\mathbb{Q})$. But this is a contradiction due to the lemma \eqref{lem:3.1}.

        \end{itemize}

\end{proof}

Now we are ready to prove the Theorem \eqref{thm:1.2}.

\begin{proof}
    From the Proposition \eqref{prop:4.1}, we have two linearly independent points namely $A_m$ and $B_m$ in $E_{m}(\mathbb{Q})$. Now from the Theorem \eqref{thm:4.1}, we can say that the cardinality of the group $E_{m}(\mathbb{Q})/ 2 E_{m}(\mathbb{Q}) $ is $2^r$, where $r $ is the rank of the group $ E_{m}(\mathbb{Q})$.  Using Proposition \eqref{prop:4.2}, we can say that $E_m(\mathbb{Q})$ has at least $4$ points which means the rank $r$ of $E_{m}(\mathbb{Q})$ is at least $2$. This concludes the proof of the Theorem.
\end{proof}

Now we give an example related to the Theorems \eqref{thm:1.1} and \eqref{thm:1.2}.
\begin{exa}
    Let us take $m=2$, $p=3$, $q=7$ and $r=11$. As these values satisfy the hypothesis of the Theorem \eqref{thm:1.1}, we can say $E_{2}(\mathbb{Q})_{tors}= \{\mathcal{O}\}$. Moreover, the rank of $E_{2}(\mathbb{Q})= 2$ that verifies Theorem \eqref{thm:1.2}. The computation was done using SAGE\cite{Sa}.
\end{exa}

\section{Concluding Remarks}
In this article, we have taken a certain family of elliptic curves $E_m$ given by the equation \eqref{eq:2.1} and shown that they are of rank at least $2$. Now using computations, it can be verified that some of these curves have rank at least $3$. So that natural question would be the following: Does there exist a sub-family of Em that has a rank at least $3$ or higher? It
would make an interesting problem.



\section*{Acknowledgement} I am grateful to the anonymous referee for his/her valuable suggestions and for helping me to improve this draft.

\end{document}